\documentclass{amsart}
\usepackage[english]{babel}
\usepackage{iftex}
\usepackage[all,cmtip]{xy}
\usepackage{hyperref}
\usepackage{lmodern}
\usepackage[initials]{amsrefs}
\usepackage{graphicx}
\usepackage[english]{babel}
\usepackage{amssymb}
\usepackage{amsthm}
\usepackage{bm}
\usepackage{enumerate}
\usepackage{stmaryrd}
\usepackage{xcolor}

\newtheorem{theorem}{Theorem}[section]
\newtheorem{lemma}[theorem]{Lemma}
\newtheorem{proposition}[theorem]{Proposition}
\newtheorem{corollary}[theorem]{Corollary}

\theoremstyle{definition}
\newtheorem{definition}[theorem]{Definition}
\newtheorem{example}[theorem]{Example}

\theoremstyle{remark}
\newtheorem{remark}[theorem]{Remark}
\newtheorem{notation}[theorem]{Notation}

\numberwithin{equation}{section}

\newcommand*{\Rmod}{R\text{-}\textbf{Mod}}
\newcommand*{\RCmod}{R\mathcal C\text{-}\textbf{Mod}}

\newcommand*{\RSmod}{R         S\text{-}\textbf{Mod}}

\newcommand*{\C}{\mathcal C} 

\newcommand*{\A}{\mathcal A}
\newcommand*{\Ab}{\textbf{Ab}}

\DeclareMathOperator{\Ker}{Ker}
\DeclareMathOperator{\Coker}{Coker}
\DeclareMathOperator{\Img}{Im}
\DeclareMathOperator{\id}{id}
\DeclareMathOperator{\Ob}{Ob}
\DeclareMathOperator{\supp}{Supp}
\DeclareMathOperator{\ind}{ind}

\DeclareMathOperator{\res}{res}
\DeclareMathOperator{\Hom}{Hom}
\DeclareMathOperator{\Mor}{Mor}
\DeclareMathOperator{\colim}{colim}
\DeclareMathOperator{\mub}{mub}

\begin{document}

\title{Generalized persistence and graded structures}

\author{Eero Hyry}
\email{eero.hyry@tuni.fi}
\address{Faculty of Information Technology and Communication Sciences,
         Tampere University,
				 Kanslerinrinne 1 (Pinni B),
         Tampere, 33100,
				 Finland
				 }

\author{Markus Klemetti}
\email{markus.o.klemetti@gmail.com}
\address{Faculty of Information Technology and Communication Sciences,
         Tampere University,
				 Kanslerinrinne 1 (Pinni B),
         Tampere, 33100,
				 Finland
				 }
\thanks{The second author was supported in part by Finnish Cultural Foundation.}


\date{February 12, 2021}


\begin{abstract}
We investigate the correspondence between generalized persistence modules and graded modules in the case the indexing set
has a monoid action. We introduce the notion of an action category over a monoid graded ring. We show that the category of additive functors from this
category to the category of Abelian groups is isomorphic to the category of modules graded over the set with a monoid action, and to the category 
of unital modules over a certain smash product. Furthermore, when the indexing set is a poset, we provide a new characterization for a generalized 
persistence module being finitely presented.\end{abstract}

\maketitle


\section{Introduction}

One of the main methods of topological data analysis is persistent homology. In the simplest case, data is encoded in an increasing nested sequence of simplicial complexes. This filtration reflects the topological and geometric structure of the data at different scales. By taking homology with coefficients in a field, one obtains the corresponding persistence module - a sequence of vector spaces and linear maps. Carlsson and Zomorodian \cite{Zomorodian}*{p.~259, Thm.~3.1~(Correspondence)} realized that one can view persistence modules as modules over a polynomial ring of one variable. The variable acts on the module as a shift. Considering filtrations indexed by $\mathbb N^n$ leads to the so called multipersistence. In~\cite{Carlsson}*{p.~78, Thm.~1}, Carlsson and Zomorodian showed that multipersistence modules now correspond to modules over a polynomial ring of $n$ variables. More generally, one can start from a filtration of a topological space indexed by a preordered set. However, the resulting generalized persistence modules do not necessarily have an immediate expression as a module over a graded ring. 

The correspondences by Carlsson and Zomorodian opened the graded perspective in topological data analysis, leading many researchers to utilize graded module theory in their investigations (see, for example,~\cite{Corbet}, \cite{Miller}, \cite{Bubenik}, \cite{Harrington}, \cite{Knudson}, \cite{Lesnick}, \cite{Scolamiero}). The most general cases of modules over a ring in this line of research are modules graded over Abelian groups with monoids as their positive cones, and modules canonically graded over cancellative monoids. 
In this article, we want to propose a new generic theoretical framework for understanding generalized persistence modules under the lens of graded algebra by considering monoid actions on preordered sets. Secondly, we want to investigate finitely presented generalized persistence modules. In particular, we will give a certain subclass of preordered sets over which finite presentation can be characterized by a suitable 'tameness' condition.

We now want to explain this in more detail. Using the language of category theory, it is convenient to define a generalized persistence module as a functor from a preordered set $P$ to the category of 
$k$-vector spaces, where $k$ is a field. In representation theory, given a commutative ring $R$ and a small category $\C$, a functor $\C\rightarrow \Rmod$ is called  an $R\C$-module. 
In this terminology, a generalized persistence module is then a $kP$-vector space. Following Mitchell~(\cite{Mitchell}), we also regard a small preadditive category $\A$ as a `ring with several objects', and an additive functor $\A \rightarrow {\Ab}$ as an $\A$-module. The $R\C$-modules may then be seen as modules over the linearization $R\C$, where $R\C$ is a preadditive category with the same objects as $\C$ and morphisms $R[\Mor_{\C}(c,d)]$, where $c,d\in \Ob \C$ (for any set $S$, we denote by $R[S]$ the free $R$-module generated by $S$).

Suppose now that $G$ is a monoid. To any $G$-act $A$, we can associate an action category $G{\smallint} A$, whose objects are the elements of $A$ and for any $a,b\in A$ the morphisms $a \rightarrow b$ are pairs $(a,g)$ where $g\in G$ with $b=ga$. It is easy to see that the category of $R(G{\smallint} A)$-modules is now equivalent to the category of $A$-graded $R[G]$-modules. Note that if the action of $G$ on $A$ is free, then simply $G{\smallint} A = A$.

Given any $G$-graded ring $S$, this leads us is to investigate the relationship between $A$-graded $S$-modules and modules over $A$ in general. We define the action category over $S$, denoted by $G{\smallint}_S A$, with objects $A$ and morphisms 
\[
s\in \bigoplus_{g\in G, \ ga=b}S_g,
\]
where $a,b\in A$. In the case $S=R[G]$, $G{\smallint}_S A$ is just the linearization of $G{\smallint} A$. Our first main result, Theorem~\ref{yleisin ekvivalenssi}, then says that the categories of $A$-graded 
$S$-modules and $G{\smallint}_S A$-modules are isomorphic. 

We can also look at the category algebra $R[G{\smallint} A]$. If $\C$ is any category, then the category algebra $R[\C]$ is defined as the free $R$-module with a basis consisting of the morphisms of $\C$, and the product of two basis elements is given by their composition, if defined, and is zero otherwise. It now turns out in Proposition~\ref{smash_isom} that the category algebra $R[G{\smallint} A]$ coincides with the smash product $R[G]\#A$, which has been much studied in ring theory (see~\cite{Nastasescu2}). This leads us to Theorem~\ref{smash_ekvivalenssi}, where we identify $G{\smallint}_S A$-modules with the category of unital $S\#A$-modules 
i.e.~the category of $S\# A$-modules $M$ with $M=(S\# A)M$. 

We then turn to consider finitely presented generalized persistence modules. Note that being finitely presented is a categorical property, so an equivalence between generalized persistence modules and graded modules preserves this property. Recall first that an 
$R\C$-module $M$ is finitely presented if there exists an exact sequence
\[
\bigoplus_{j\in J} R[\Mor_{\C}(d_j,-)] \rightarrow \bigoplus_{i\in I} R[\Mor_{\C}(c_i,-)] \rightarrow M\rightarrow 0,
\]
where $I$ and $J$ are finite sets, and $c_i,d_j\in \C$ for all $i\in I$, $j\in J$.
We will look at posets $\C$ which are weakly bounded from above and mub-complete. By `weakly bounded' we mean that every finite subset $S \subseteq \C$ has a finite number of minimal upper bounds in $\C$, whereas $\C$ is mub-complete if given a finite non-empty subset $S \subseteq \C$ and an upper bound $c$ of $S$, there exists a minimal upper bound $s$ of $S$ such that $s \le c$. 
In our Theorem~\ref{fin_pres}, we characterize finitely presented generalized persistence modules in this situation. More precisely, we can show that an $R\C$-module $M$ is finitely presented if and only if the 
$R$-modules $M(c)$ are finitely presented for all $c\in \C$, and $M$ is $S$-determined for some finite set $S\subseteq\C$. 

Given $S\subseteq\C$, we call an $R\C$-module $M$ $S$-determined if  $\supp M \subseteq{\uparrow} S$ and the implication 
\[
S\cap {\downarrow} c =  S \cap {\downarrow} d \ \Rightarrow \ \text{the morphism}\ M(c\leq d) \text{ is an isomorphism}
\]
holds for every $c\le d$ in $\C$. Here $\supp(M):=\{c\in \C\mid M(c)\neq 0\}$ denotes the support of $M$,
and for any  $T\subseteq\C$, we use the usual notations 
\[
\hbox{${\uparrow} T := \{c\in \C \mid t\leq c \text{ for some }t\in T\}$ and
${\downarrow} T := \{c\in \C \mid c\leq t \text{ for some }t\in T\}$}
\]
for the upset generated and the downset cogenerated by $T$, respectively. 
Our intuition for this definition comes from topological data analysis, where one tracks how the elements of 
each $M(c)$ evolve in the morphisms $M(c\leq c')$ ($c,c'\in \C)$. One says that an element $m\in M(c)$ is born at $c$ if it is not in the image of any morphism $M(c'\le c)$, where $c'<c$, and dies at $c''$ 
if $M(c\le c'')(m)=0$ and $M(c\leq c')(m)\neq 0$ for all $c\leq c'<c''$. Suppose that there exists a set $S$ such that all births and deaths occur inside $S$. The condition $S\cap {\downarrow} c =  S \cap {\downarrow} d$ then implies that looking down from both $c$ and $d$, we see the same deaths and births. In particular, the morphism $M(c\leq d)$ must be an isomorphism.

Our proof for Theorem~\ref{fin_pres} starts from the fact that an $R\C$-module $M$ is finitely presented if and only if the $R$-modules $M(c)$ are finitely presented for all $c\in \C$ and $M$ is $S$-presented for some finite subset $S\subseteq \C$. Here $S$-presented means the existence of a set $S\subseteq \C$ and an exact sequence of the type 
\[
\bigoplus_{s\in S}B_s[\Mor_\C(s,-)]  \rightarrow  \bigoplus_{s\in S}A_s[\Mor_\C(s,-)] \rightarrow  M  \rightarrow  0,
\]
where $A_s$ and $B_s$ are $R$-modules for all $s\in S$. It is easily seen that if $M$ is $S$-presented, then $M$ is $S$-determined.  We denote the set of minimal upper bounds of non-empty subsets of a finite set $S\subseteq \C$ by
\[
\hat{S} := \bigcup_{\emptyset\neq S'\subseteq S} \mub_\C(S').
\]
In Corollary~\ref{tuplahattu} we now make the crucial observation that $M$ is $\hat{\hat S}$-presented if $S\subseteq\C$ is a finite set such that $M$ is $S$-determined.

As a useful tool we introduce the sets of births and and deaths relative to $S$ by
\[
B_S(M):=\{c\in \C\mid \colim_{s<c, \ s\in S} M(s) \rightarrow M(c) \ \text{is a non-epimorphism}\}
\]
and
\[
D_S(M):=\{c\in \C\mid \colim_{s<c, \ s\in S} M(s) \rightarrow M(c) \ \text{is a non-monomorphism}\}.
\]   
An $R\C$-module $M$ is known to be $S$-presented if and only if the natural morphism $\ind_S\res_S M\rightarrow M$ is an isomorphism. Here
$\res_S$ denotes the restriction functor  from the category of $R\C$-modules to the category of $RS$-modules, and $\ind_S$ its left adjoint, 
the induction functor, in the opposite direction. Note the pointwise formula 
\[
(\ind_S\res_S M)(c)= \colim_{s\leq c, \ s\in S} M(s)
\]
for all $c\in \C$. 
We observe in Proposition~\ref{FSP_apu} that the module $M$ is $S$-presented if and only if $B_S(M)\cup D_S(M) \subseteq S$. Interestingly, if $S$ is Artinian, then
$B_S(M)\cup D_S(M)$ is the minimal subset $T\subseteq S$ such that $M$ is $T$-presented (see Proposition~\ref{esitys_minimi_osa2}). 
Suppose that $\C=\mathbb Z^n$, $R=k$ is a field and 
\[
0\rightarrow L \rightarrow N\stackrel{f}{\rightarrow} M\rightarrow 0,
\]
is an exact sequence, where $N$ is a free module and $f$ a minimal epimorphism. In this case our Theorem~\ref{verho} says that $D_S(M)=B_S(L)$ confirming the intuition that deaths should correspond to
`relations'.

This work unifies several earlier results. In the context of topological data analysis, monoid actions have been considered by Bubenik et al.~in their article~\cite{Bubenik2}, where they looked at the 
action on any preordered set given by the monoid of its translations. In the article~\cite{DeSilva} of de Silva et al., an indexing category with an additional structure of a $[0,\infty)$-action is called a category with a coherent flow. Recently, Bubenik and Milicevic considered modules graded over Abelian groups with monoids as their positive cones~(\cite{Bubenik}). We have in particular been motivated by the article~\cite{Corbet} of Corbet and Kerber, who generalized the result of Carlsson and Zomorodian to the case where the indexing set is a so called good monoid.  We point out that if $G$ is a monoid, then $RG$-modules of finitely presented type of Corbet and Kerber (\cite{Corbet}*{p.~19, Def.~15}) are the same thing as finitely presented $RG$-modules. The set $\hat S$ is a framing set
in the sense of~\cite{Corbet}*{p.~19, Def.~15}. If the set $S$ happens to be a closed interval in $\mathbb Z^n$, then, up to translation, `$S$-determined' means the same as `positively $a$-determined' for some $a\in \mathbb N^n$ as defined by Miller in~\cite{Miller3}*{p.~186, Def.~2.1}. Note that this is a special instance of `finite encoding' investigated in \cite{Miller}. Our sets of births and and deaths relative to $S$ are related to the invariants $\xi_0$ and $\xi_1$ studied by Carlsson and Zomorodian in~\cite{Carlsson}, and also by Knudson in~\cite{Knudson}. For a finitely generated $\mathbb Z^n$-graded $k[X_1,\ldots,X_n]$-module $M$, the invariants $\xi_0(M)$ and $\xi_1(M)$ are multisets indicating the degrees of minimal generators and minimal relations of $M$ equipped with the multiplicities they occur. The underlying sets of $\xi_0(M)$ and $\xi_1(M)$ are now $B_S(M)$ and $D_S(M)$. 

We assume that the reader is familiar with the basic notions of category theory, see for example~\cite{MacLane}. For more details on $R\C$-modules, we refer to~\cite{Luck} and~\cite{TomDieck}.


\section{Modules over a monoid act}

\subsection{Monoid action}

In this section we recall some basic properties of monoid actions. Let $G$ be a monoid and let $A$ be a set. If there exists an operation $\cdot: G\times A\rightarrow A$ such that $(gh)a = g(ha)$ and $1_G \cdot a = a$ for all $g,h\in G$ and $a\in A$, we say that $A$ is a \emph{(left) $G$-act}. We then get a preorder on $A$ by setting $a\leq b$ if $b=ga$ for some $g \in G$. The action naturally gives rise to two categories having $A$ as the set of objects. 

First, we have a small thin category $A$, where for all $a,b\in A$ there exists a unique morphism $a\rightarrow b$ if $a\leq b$ in the preorder. By abuse of notation we write $a\leq b$ for this morphism. Recall that in general a category is \emph{thin} if there exists at most one morphism between any two objects. 

Secondly, there is an \emph{action category} $G {\smallint} A$, where morphisms $a\rightarrow b$ are pairs $(a,g)$ such that $b=ga$ for some $g\in G$. If there is no possibility of confusion, we sometimes denote the morphism $(a,g)$ by 
$g$. Composition of morphisms in $G{\smallint} A$ is defined by the multiplication of $G$:
\[
(ga,h)\circ (a,g) = (a,hg).
\]

There is an obvious functor $G{\smallint} A\rightarrow A$ where
\[
a\mapsto a \quad \text{and} \quad (a,g)\mapsto (a\leq ga).
\]
This functor is an isomorphism if and only if the $G$-action on $A$ is \emph{free}, i.e.~for all $g,h\in G$,
\[
ga=ha \ \text{for some}\ a\in A \ \Rightarrow \ g=h.
\]

\begin{remark}\label{elements}
We often consider the monoid $G$ itself as a $G$-act, so it gives rise to a thin category $G$ and the action category $G\smallint G$. Sometimes, the monoid $G$ is viewed as a category $BG$ with a single object. Then a monoid act could be defined as a functor $F:BG \rightarrow \textbf{Set}$. The corresponding action category coincides with the category of elements of $F$ (see \cite{Riehl2}*{p.~66, Ex.~2.4.10}).
\end{remark}

\begin{example}\label{ryhmatoiminta} An abelian group $G$ is called \emph{preordered} if it is equipped with a preorder $\leq$ such that $g\leq g'$ implies $g+h\leq g'+h$ for all
$g,g',h\in G$. If $G_+=\{g\in G\mid g\geq 0\}$ is its \emph{positive cone}, then $g\leq g'$ is equivalent to $g'-g\in G_+$.
The action of the monoid $G_+$ on $G$ is free, and we may identify the action category $G_+ {\smallint} G$ with $G$. 
\end{example}

A \emph{translation} on a preordered set $P$ is an order-preserving function $F:P\rightarrow P$ that satisfies the condition $p\leq F(p)$ for all $p\in P$. The translations of $P$ form a monoid $\text{Trans}(P)$ with composition as the operation. 

Let $A$ be a $G$-act. The action by an element $g\in G$ now determines a translation on $A$ if and only if $a\leq b$ implies $ga\leq gb$.
If this implication holds for all $g\in G$, then we say that $A$ is an \emph{order-preserving $G$-act}. Note that any $G$-act $A$ is order-preserving if $G$ is commutative. For an order-preserving $G$-act $A$, we get a monoid homomorphism $\varphi$ from $G$ into the monoid of translations $\text{Trans}(A)$. This induces a monoid embedding $\hat{\varphi}:G/\Ker \varphi \rightarrow \text{Trans}(A)$, where $\Ker\varphi$ is the congruence relation defined by
\[
(g,h)\in \Ker \varphi \ \Leftrightarrow \ ga=ha \ \text{for all} \ a\in A.
\]
In particular, $\varphi$ is an embedding if and only if the $G$-action on $A$ is \emph{faithful}: for all $g,h\in G$, 
$ga=ha$ for all $a\in A$ implies that $g=h$. 

We next give a slight generalization of~\cite{Garcia}*{p.~4, Thm.~2.2}.

\begin{proposition}\label{indeksit}
For any preordered set $P$, there exists a monoid $G$ and a $G$-act $A$ such that $P$ and $A$ are isomorphic as thin categories.
\end{proposition}

\begin{proof}
We present the proof here for the convenience of the reader. Let $G$ denote the submonoid of the monoid of all functions $P\rightarrow P$ consisting 
of the functions $g:P\rightarrow P$ for which $a\leq_P g(a)$ for all $a\in P$.  Define the $G$-action on $A:=P$ by setting $g\cdot a = g(a)$ for all
$g\in G$ and $a\in A$. Then $A$ is a $G$-act. It remains to show that $a\leq_P b$ if and only if $a\leq_A b$.

Assume first that $a\leq_A b$. By definition, there exists an element $g\in G$ such that $b=ga$. But this means that 
$a\leq_P g(a)=ga=b$.
Conversely, if $a\leq_P b$, we define a function $g:P\rightarrow P$ by setting
\[
g(p)=\left\{\begin{aligned}
&b, \ \text{if} \ p=a;\\
&p, \ \text{otherwise}.
\end{aligned}\right.
\]
We then immediately see that $g\in G$ and $ga=g(a)=b$, so that $a\leq_A b$
\end{proof}

\begin{remark}\label{uniq}
The monoid $G$ in Proposition~\ref{indeksit} does not need to be unique: For example, the one element set $P=\{*\}$ has the trivial monoid action for any monoid $G$. Note that if $G$ is the monoid of the proof of Proposition~\ref{indeksit}, then the order-preserving elements in $G$ are exactly the translations on $A$, so that $\text{Trans}(A)\subseteq G$. 
\end{remark}


\subsection{Action categories over a graded ring}

Theorem~\ref{yleisin ekvivalenssi} will generalize the equivalence of the correspondence theorem of Carlsson and Zomorodian \cite{Zomorodian}*{p.~259, Thm.~3.1} mentioned in the Introduction. Moreover, it generalizes the multi-parameter version of the theorem by Carlsson and Zomorodian~(\cite{Carlsson}*{p.~78, Thm.~1}) as well as the generalization given by Corbet and Kerber~(\cite{Corbet}*{p.~18, Lemma~14}). For a discussion on related finiteness conditions, see~\cite{Corbet}*{p.~3} and Remark~\ref{scola_corbet}.

We begin by defining a certain preadditive category.

\begin{definition}
Let $A$ be a $G$-act, and let $S:=\oplus_{g\in G} S_g$ be a $G$-graded ring. The \emph{action category over $S$}, denoted $G{\smallint}_S A$, is the category with the set $A$ as objects, and morphisms $(a,s)\colon  a\rightarrow b$, where $a,b\in A$ and
\[
s\in \bigoplus_{g\in G, \ ga=b}S_g.
\]
Composition for morphisms $(a,s)\colon a\rightarrow ga$ and $(ga,t)\colon  ga\rightarrow hga$ is defined by
\[
(ga,t)\circ(a,s)  = (a,ts).
\]
\end{definition}

\begin{remark}\label{huom4}
Keeping a close eye on the domains, we may write $s:=(a,s)$. With this notation, composition is just the multiplication in $S$.
\end{remark}

\begin{example}\label{lin_gen}
Let $A$ be a $G$-act. If $R$ is a commutative ring, then the action category $G{\smallint}_{R[G]} A$ over the monoid ring $R[G]$ coincides with 
the linearized action category $R(G{\smallint} A)$. Indeed, by definition $\Ob R(G{\smallint} A) = A$, and
\[
\Hom_{R(G {\smallint} A)}(a,b)=R[\{(a,g) \mid  g\in G \text{ and }ga=b\}].
\]
for all $a,b\in A$.
\end{example}

\begin{example}\label{companion}
If $G$ is an Abelian group and $S:=\bigoplus_{g\in G}S_g$ is a $G$-graded ring, the category $G{\smallint}_S G$ is called in~\cite{DellAmbrogio}*{p.~358, Def.~2.1} a \emph{companion category}. In this case, we may identify $\Hom_{G{\smallint}_S G}(g,h)$ with $S_{h-g}$.
\end{example}

Let $A$ be a $G$-act. Let $S:=\bigoplus_{g\in G}S_g$ be a $G$-graded ring. Recall that a (left) $S$-module $M$ is \emph{$A$-graded}, if
\begin{enumerate}
\item[1)] $M = \bigoplus_{a\in A} M_a$, where $M_a$ is an Abelian group for all $a\in A$;
\item[2)] $S_g M_a \subseteq M_{ga}$ for all $g\in G$ and $a\in A$.
\end{enumerate}

Preparing for Theorem~\ref{yleisin ekvivalenssi}, we will now define two functors, $\Phi$ and $\Psi$, that connect $A$-graded $S$-modules to $(G{\smallint}_S A)$-modules. 

Let $M$ be a $G{\smallint}_S A$-module. By setting $sm = M(s)(m)$ for all $g\in G$, $s\in S_g$ and $m\in M(a)$, we can define an $A$-graded $S$-module
\[
\Phi M := \bigoplus_{a\in A} M(a).
\]
A morphism $f\colon M\rightarrow N$ of $G\smallint_S A$-modules consists of homomorphisms of Abelian groups $f_a\colon  M(a)\rightarrow N(a)$ with commutative diagrams
\[
\xymatrix{
M(a)\ar[d]_{f_a}\ar[r]^{M(s)}&  M(ga) \ar[d]^{f_{ga}}\\
N(a) \ar[r]_{N(s)}& N(ga)}
\] 
for all $a\in A$, $g\in G$ and $s\in S_g$. These homomorphisms and diagrams obviously give rise to a homomorphism $\Phi f\colon \Phi M\rightarrow \Phi N$ of $A$-graded $S$-modules with $(\Phi f)_a=f_a$
for all $a\in A$.

Next, let $Q$ be an $A$-graded $S$-module. We set $(\Psi Q)(a)=Q_a$ for all $a\in A$. If $(a,s)\colon a\rightarrow ga$ is a morphism, where $a\in A$, $g\in G$ and $s\in S_g$, we can define a homomorphism 
\[
(\Psi Q)((a,s))\colon (\Psi Q)(a) \rightarrow (\Psi Q)(ga)
\]
 by setting $(\Psi Q)((a,s))(q) = s\cdot q$ for all $q\in Q_a$. It is clear that $\Psi Q$ is an additive functor $G\smallint_S A \rightarrow {\Ab}$, i.e.,~a $G\smallint_S A$-module. Moreover, if $h\colon Q\rightarrow P$ is a homomorphism of $A$-graded $S$-modules, we have a morphism of $G\smallint_S A$-modules $\Psi h\colon ~\Psi Q\rightarrow \Psi P$ given 
by $(\Psi h)_a = h_a$ for all $a\in A$.

\medskip

We are now ready to state

\begin{theorem}\label{yleisin ekvivalenssi}
Let $A$ be a $G$-act, and let $S:=\oplus_{g\in G} S_g$ be a $G$-graded ring. The above functors $\Phi$ and $\Psi$ give an isomorphism of categories
\[
(G{\smallint}_S A)\text{-}\textbf{Mod} \cong A\text{-}\text{gr } S\text{-}\textbf{Mod}.
\]
\end{theorem}

\begin{proof}
It remains to prove that $\Phi\circ \Psi=\id$ and $\Psi\circ \Phi=\id$, which is straightforward.
\end{proof}

Combining this theorem with Example~\ref{lin_gen} gives 

\begin{corollary}\label{lineaari ekvivalenssi}
Let $A$ be a $G$-act, and let $R$ be a commutative ring. There is an isomorphism of categories
\[
R(G{\smallint} A)\text{-}\textbf{Mod} \cong A\text{-}\text{gr } R[G]\text{-}\textbf{Mod}.
\]
In particular, if the $G$-action on $A$ is free, we obtain an isomorphism
\[
RA\text{-}\textbf{Mod} \cong A\text{-gr } R[G]\text{-}\textbf{Mod}.
\]
\end{corollary}

\begin{example}\label{triviaali ekvivalenssi}
If $A=\{e\}$ is a one object set, Theorem~\ref{yleisin ekvivalenssi} gives us an isomorphism $G{\smallint}_S\{e\}\text{-}\textbf{Mod}\cong S\text{-}\textbf{Mod}$. In the case $S=R[G]$, where $R$ is a commutative ring, this means that $RG\text{-}\textbf{Mod}\cong R[G]\text{-}\textbf{Mod}$, where $RG$ is the linearization of the $1$-object category $G$.
\end{example}

\begin{example}Let $G$ be a preordered Abelian group with the positive cone $G_+$ (see Example~\ref{ryhmatoiminta}). If $R$ is a commutative ring, then by Corollary~\ref{lineaari ekvivalenssi} the 
categories $RG\text{-}\textbf{Mod}$ and  $G\text{-}\text{gr }R[G_+]\text{-}\textbf{Mod}$ are isomorphic.
\end{example}


\subsection{Category algebras and smash products}

Let $\C$ be a small category, and let $R$ be a commutative ring. A \emph{category algebra} $R[\C]$ is the free $R$-module with the basis consisting of the elements $e_u$, where $u\colon c\rightarrow d$ is a morphism in $\C$, and with multiplication defined by
\[
e_v\cdot e_u = 
\left\{
\begin{aligned}
&e_{vu}, \ \text{if} \ c'=d;\\
&0, \ \text{otherwise}
\end{aligned}
\right.
\] 
for morphisms $u\colon c\rightarrow d$ and $v\colon c'\rightarrow d'$ in $\C$. Equipped with this product, $R[\C]$ becomes a ring that has a unit if $\C$ is finite. 

Let $A$ be a $G$-act and $S$ a $G$-graded ring. We recall (see \cite{Nastasescu}*{p.~390}) that a \emph{smash product} $S\#A$ is the free (left) $S$-module with the basis $\{p_a\mid a\in A\}$, and with multiplication defined by the bilinear extension of
\[
(s_g p_a)(t_h p_b) = \left\{
\begin{aligned}
&(s_g t_h)p_b, \ \text{if} \ hb=a;\\
&0, \ \text{otherwise}
\end{aligned}
\right.
\]
where $g,h\in G$, $s_g\in S_g$, $t_h\in S_h$ and $a,b\in A$. Equipped with this multiplication, $S\# A$ is a non-unital ring, i.e.~a ring possibly without identity. However, $S\#A$ has \emph{local units}. This means that every finite subset of $S\# A$ is contained in a subring of the form $w(S\# A)w$, where $w$ is an idempotent of $S\# A$. More precisely, let $T:=\{t_1,\ldots,t_n\}$ be a finite subset of $S\# A$. We may assume that $t_i= s_i p_{a_i}$, where $g_i\in G$, $a_i\in A$ and $s_i\in S_{g_i}$ for all $i\in \{1,\ldots,n\}$. We denote
\[
B := \{ a \in A \mid a=a_i \text{ or } a=g_ia_i \text{ for some } i\in \{1,\ldots,n\}\}
\]
and $w:=\sum_{a\in B} p_a$. It is now straightforward to see that $w$ is idempotent and $wt_iw = wt_i = t_i$ for all $i\in \{1,\ldots,n\}$. 

Let $R'$ be a non-unital ring. An $R'$-module $M$ is \emph{unital} if it satisfies the condition $M=R'M$.

The next proposition and its proof are inspired by~\cite{Beattie}*{p.~221, Cor.~2.4}.

\begin{proposition}\label{beattie221}
Let $M$ be an $S\#A$-module. Then $M$ is unital if and only if for every finite subset $N\subseteq M$ there exists a finite subset $B\subseteq A$ such that $wn=n$ for all $n\in N$, where $w:=\sum_{a\in B} p_a$.
\end{proposition}

\begin{proof}
Assume first that $M$ is unital. Let $N:=\{n_1,\ldots,n_p\}\subseteq M$ be a finite set. Now, for all $i\in\{1,\ldots,p\}$, the element $n_i$ may be written as
\[
n_i=\sum_{j=1}^q s_{i,j}n_{i,j},
\]
where $s_{i,j}\in S\# A$ and $n_{i,j}\in M$ for all $j\in\{1,\ldots,q\}$. This gives us a finite set
\[
T=\{s_{i,j}\mid i\in\{1,\ldots,p\}, \ j\in\{1,\ldots,q\}\}\subseteq S\# A.
\]
As stated above, we then have a finite subset $B\subseteq A$ such that $w=ws$ for all $s\in T$, where $w:=\sum_{a\in B} p_a$. Thus for all $i\in\{1,\ldots,p\}$,
\[
wn_i=w(\sum_{i=1}^q s_{i,j}n_{i,j}) = \sum_{i=1}^q(ws_{i,j})n_{i,j} = \sum_{i=1}^q s_{i,j}n_{i,j} =n_i.
\]

Conversely, suppose that for every finite subset $N\subseteq M$ there exists a finite subset $B\subseteq A$ such that $wn=n$ for all $n\in N$, where $w:=\sum_{a\in B} p_a$. Taking $N=\{m\}$ for $m\in M$, we get $m=wm\in S\# A$.
\end{proof}

\begin{proposition}\label{smash_isom}
Let $R$ be a commutative ring, $G$ a monoid,  and $A$ a $G$-act. There exists an isomorphism of non-unital rings
\[
\varphi\colon  R[G{\smallint} A] \rightarrow R[G]\# A
\]
defined by $e_{(a,g)} \mapsto e_g p_a$ for all $a\in A$ and $g\in G$.
\end{proposition}

\begin{proof}
It is easy to see that $\varphi$ is an isomorphism of $R$-modules. It is also a ring homomorphism, since for all $a,b\in A$ and $g,h\in G$,
\begin{align*}
\varphi(e_{(b,h)}e_{(a,g)})&=\left\{\begin{aligned}  &\varphi(e_{(a,hg)}), \text{ if } b=ga;\\ &0, \ \text{else} \end{aligned}\right.\\
&=\left\{\begin{aligned}  &e_{hg} p_a, \text{ if } b=ga;\\ &0, \ \text{else} \end{aligned}\right.\\
&=(e_hp_b)(e_gp_a)\\
&=\varphi(e_{(b,h)})\varphi(e_{(a,g)}).
\end{align*}
\end{proof}

\begin{proposition}\label{smash_unital}
Let $M$ be an $S\# A$-module. Then $M=\bigoplus_{a\in A}p_a M$ if and only if $M$ is unital.
\end{proposition}

\begin{proof}
Assume first that $M=\bigoplus_{a\in A}p_a M$. Let $N:=\{n_1,\ldots,n_p\}\subseteq M$. Since for all $i\in\{1,\ldots,p\}$, the element $n_i$ may be written as
\[
n_i=\sum_{j=1}^q p_{a_{i,j}}n_{i,j},
\]
where $a_{i,j}\in A$ and $n_{i,j}\in M$ for all $j\in\{1,\ldots,q\}$, there exists a finite subset 
\[
B:=\{a_{i,j} \mid i\in\{1,\ldots,p\}, \ j\in\{1,\ldots,q\}\}
\]
of $A$. Let $w:= \sum_{a\in B} p_{a}$. Then
\[
wn_i = w\left(\sum_{j=1}^q p_{a_{i,j}}n_{i,j}\right) = \sum_{j=1}^q wp_{a_{i,j}}n_{i,j}= n_i,
\]
so $M$ is unital by Proposition~\ref{beattie221}.

Assume next that $M$ is unital. Let $m\in M$. By Proposition~\ref{beattie221}, we may write $m=wm$ for some $w=\sum_{a\in B}p_a$, where $B\subseteq A$ is finite. Thus
\[
m=(\sum_{a\in B}p_{a})m = \sum_{a\in B} p_{a}m,
\]
so that $M=\sum_{a\in A} p_{a} M$. Furthermore, since the elements $p_{a}$ are orthogonal, the sum is direct.  
\end{proof}

Let us denote by $S\# A$-\textbf{Mod} the category of unital $S\# A$-modules. We will now define two functors, $\Gamma$ and $\Lambda$, that connect unital $(S\# A)$-modules to $(G{\smallint}_S A)$-modules. Let $M$ be a $G{\smallint}_S A$-module. Set
\[
\Gamma M := \bigoplus_{a\in A}M(a).
\]
It is not difficult to check that by setting $(sp_a)m = M((a,s))(m_a)$ for all $g\in G$, $s\in S_g$, $a\in A$ and $m:=\sum_{b\in A}m_b\in \Gamma M$, $\Gamma M$ becomes an $S\# A$-module.
To show unitality, notice that $p_a(\Gamma M) = M(a)$ for all $a\in A$, which implies that $\Gamma M=\bigoplus_{a\in A}p_a(\Gamma M)$. Thus $\Gamma M$ is unital by Proposition~\ref{smash_unital}. If $f: M\rightarrow N$ is a morphism of $G{\smallint}_S A$-modules, we can define a homomorphism $\Gamma f: \Gamma M \rightarrow \Gamma N$ of $(S\# A)$-modules by setting 
\[(\Gamma f)(m)=\sum_{a\in A}f_a(m_a)\]
for all $m=\sum_{a\in A}m_a\in \Gamma M$.

Next, let $Q$ be a unital $S\# A$-module. We define a $G{\smallint}_S A$-module $\Lambda Q$ by first setting $(\Lambda Q)(a)=p_aQ$ for all $a\in A$. Let $a\in A$ and $g\in G$, $s\in S_g$. Given a 
morphism $(a,s)\colon a\rightarrow ga$, we then have a homomorphism of Abelian groups
\[(\Lambda Q)((a,s))\colon (\Lambda Q)(a)\rightarrow (\Lambda Q)(ga), q\mapsto (sp_a)q.\]
Finally, for a homomorphism $h\colon Q\rightarrow P$ of $S\# A$-modules, there is a morphism of $G{\smallint}_S A$-modules $\Lambda h\colon \Lambda Q\rightarrow \Lambda P$ with $(\Lambda h)_a(q)=h(q)$ for all $a \in A$ and $q\in (\Lambda Q)(a)$.

\begin{theorem}\label{smash_ekvivalenssi}
Let $A$ be a $G$-act, and let $S:=\oplus_{g\in G} S_g$ be a $G$-graded ring. The functors $\Gamma$ and $\Lambda$ give an isomorphism of categories
\[
(G{\smallint}_S A)\text{-}\textbf{Mod} \cong S\# A \text{-}\textbf{Mod}.
\]
\end{theorem}

\begin{proof}
We need to show that $\Gamma \Lambda =\id$ and $\Lambda \Gamma =\id$.

Let $Q$ be a unital $S\# A$-module. By Proposition~\ref{smash_unital} we then have
\[
(\Gamma \Lambda)Q = \bigoplus_{a\in A} (\Lambda Q)(a) = \bigoplus_{a\in A} p_aQ = Q.   
\]
Moreover, the $S\# A$-module structures of $Q$ and $(\Gamma\Lambda) Q$ are the same. Indeed, writing $*$ for the multiplication by $S\#A$ on $(\Gamma\Lambda) Q$, we get
\[
(sp_a)*q = (\Lambda Q)((a,s))(p_a q_a) = (sp_a) (p_a q_a) = (sp_a) q
\]
for all $a\in A$, $g\in G$, $s\in S_g$ and $q:=\sum_{a\in A}p_a q_a\in Q$. 

On the other hand, let $M$ be a $G{\smallint}_S A$-module. For an object $a\in A$,
\[
((\Lambda \Gamma)M)(a)= p_a(\Gamma M) = M(a).
\]
Furthermore, if $(a,s)\colon a\rightarrow ga$ is a morphism in $G{\smallint}_S A$, then
\[
((\Lambda \Gamma)M)((a,s))(m)=(sp_a)m = M((a,s))(m)
\]
for all $m\in M(a)$, so that $((\Lambda \Gamma)M)((a,s))= M((a,s))$. 
\end{proof}

\begin{corollary}\label{smash_lin_ekvivalenssi}
Let $A$ be a $G$-act, and let $R$ be a commutative ring. There exists an isomorphism of categories between the categories of $R(G{\smallint} A)$-modules and unital $R[G{\smallint} A]$-modules.
\end{corollary}

\begin{proof}
This follows from Proposition~\ref{smash_isom} and Theorem~\ref{smash_ekvivalenssi}.
\end{proof}


\section{Finitely presented \texorpdfstring{$R\C$}{RC}-modules}

We will assume in the following that $\C$ is a small category and $R$ a commutative ring.  Recall first that an $R\C$-module $M$ is 

\begin{itemize}
\item \emph{finitely generated} if there exists an epimorphism
\[
\bigoplus_{i\in I} R[\Mor_{\C}(c_i,-)] \rightarrow M
\]
where $I$ is a finite set, and $c_i\in \C$ for all $i\in I$;
\item \emph{finitely presented} if there exists an exact sequence
\[
\bigoplus_{j\in J} R[\Mor_{\C}(d_j,-)] \rightarrow \bigoplus_{i\in I} R[\Mor_{\C}(c_i,-)] \rightarrow M\rightarrow 0,
\]
where $I$ and $J$ are finite sets, and $c_i, d_j\in \C$ for all $i\in I$ and $j\in J$.
\end{itemize}
For more details on finitely generated and finitely presented objects in an Abelian category, we refer the reader to~\cite{Popescu}*{Ch.~3.5}.


\subsection{\texorpdfstring{$S$}{S}-presented and \texorpdfstring{$S$}{S}-generated \texorpdfstring{$R\C$}{RC}-modules}

Let $S\subseteq \C$ be a full subcategory. The notions of $S$-generated and $S$-presented modules will play an important role in the rest of this article. Before going into details, we will recall some facts about the restriction and induction functors along the inclusion $i\colon S\subseteq \C$. 

The restriction $\res_S\colon \RCmod \rightarrow \RSmod$ is defined 
by precomposition with $i$, and the induction $\ind_S\colon \RSmod \rightarrow \RCmod$ is its left Kan extension along $i$. The induction is the left adjoint of the restriction. 
Note, in particular, that 
it thus commutes with colimits. The counit of this adjunction gives  us for every $R\C$-module $M$ the  \emph{canonical morphism} \[\mu_M\colon  \ind_S\res_S M\rightarrow M,\] 
which we will use frequently.

More explicitly, for any $R\C$-module $M$ and $RS$-module $N$, we have the pointwise formulas 
\[
(\res_S M)(s) = M(s) \quad \text{and}\quad (\ind_S N)(c)  = \colim_{(t,u)\in (i/c)}N(t)
\]
for all $s\in S$ and $c\in \C$. Here $(i/c)$ denotes the slice category. Its objects are pairs $(s,u)$, where $s\in S$ and $u\colon s\rightarrow c$ is a morphism in $\C$. For $(s,u),(t,v)\in \Ob
(i/c)$, a morphism $(s,u) \rightarrow (t,v)$ is a morphism $\alpha\colon s\rightarrow t$ in $S$ with $v\alpha=u$.
If $\C$ is a poset, the latter formula yields
\[
(\ind_S N)(c) = \colim_{t\in S, \ t\leq c}N(t).
\] 

Let $A$ be an $R$-module and $c\in \C$. We define an $R\C$-module
\[
A[\Mor_\C(c,-)]:= A \otimes_R R[\Mor_\C(c,-)]
\]
by taking a pointwise tensor product. We note that the functor $\Rmod\rightarrow \RCmod$ that sends $A$ to $A[\Mor_\C(c,-)]$ is right exact for all $c\in\C$.

\begin{proposition}\label{induktio_apu}
Let $S\subseteq \C$ be a full subcategory, $A$ an $R$-module, and $s\in \C$. Then 
\[
\ind_S\res_S A[\Mor_\C(s,-)]\cong A[\Mor_\C(s,-)].
\]
\end{proposition}

\begin{proof}
By Yoneda's lemma and the aforementioned adjunction, we have the following isomorphisms:
\begin{align*}
\Hom_{R\C}(R[\Mor_\C(s,-)], M) &\cong M(s)\\
&\cong \Hom_{RS}(R[\Mor_S(s,-)],\res_S M)\\
&\cong \Hom_{R\C}(\ind_S R[\Mor_S(s,-)],M).
\end{align*}
This shows us that $\ind_S \res_S R[\Mor_\C(s,-)]\cong R[\Mor_\C(s,-)]$. In particular
\[
\colim_{(t,u)\in (i/d)}R[\Mor_\C(s,t)] \cong R[\Mor_\C(s,d)]  
\]
for $d\in \C$. Since tensoring commutes with colimits, we see that for all $d\in \C$,
\begin{align*}
(\ind_S\res_S A[\Mor_\C(s,-)])(d)&= \colim_{(t,u)\in (i/d)}A[\Mor_\C(s,t)] \\
&\cong A\otimes_R \colim_{(t,u)\in (i/d)}R[\Mor_\C(s,t)] \\
&\cong A\otimes_R R[\Mor_\C(s,d)] \\
&\cong A[\Mor_\C(s,d)].
\end{align*}
Therefore $\ind_S \res_S A[\Mor_\C(s,-)] \cong A[\Mor_\C(s,-)]$ as wanted. 
\end{proof}

An $R\C$-module $M$ is said to be \emph{$S$-generated} if the natural morphism 
\[
\bigoplus_{s\in S}M(s)[\Mor_\C(s,-)]\rightarrow M
\]
is an epimorphism. Since this morphism factors through the canonical morphism $\mu_M\colon  \ind_S\res_S M\rightarrow~M$, we see that $M$ is $S$-generated if and only if $\mu_M$ is an epimorphism.

\begin{proposition}\label{s-esitetyt}
Let $S\subseteq \C$ be a full subcategory. Assume that $M$ is an $S$-generated $R\C$-module, so that we have an exact sequence of $R\C$-modules
\[
0\rightarrow K \rightarrow \bigoplus_{s\in S}M(s)[\Mor_\C(s,-)]\rightarrow M \rightarrow 0.
\]
Then the following are equivalent:
\begin{itemize}
\item[1)] The canonical morphism $\mu_M\colon  \ind_S\res_S M\rightarrow M$ is an isomorphism;
\item[2)] If there exists an exact sequence of $R\C$-modules
\[
0\rightarrow L \rightarrow N\rightarrow M \rightarrow 0,
\]
where $N$ is $S$-generated, then $L$ is $S$-generated;
\item[3)] $K$ is $S$-generated;
\item[4)] The sequence 
\[
\bigoplus_{s\in S}K(s)[\Mor_\C(s,-)] \rightarrow \bigoplus_{s\in S}M(s)[\Mor_\C(s,-)]\rightarrow M \rightarrow 0
\]
is exact;
\item[5)] For each $s\in S$, there exist $R$-modules $A_s$ and $B_s$ such that the sequence
\[
\bigoplus_{s\in S}B_s[\Mor_\C(s,-)] \rightarrow \bigoplus_{s\in S}A_s[\Mor_\C(s,-)]\rightarrow M \rightarrow 0
\] 
is exact.
\end{itemize}
When these equivalent conditions hold, we say that $M$ is \emph{$S$-presented}.
\end{proposition}

\begin{proof}
We will show that 1) $\Rightarrow$ 2) $\Rightarrow$ 3) $\Rightarrow$ 4) $\Rightarrow$ 5) $\Rightarrow$ 1). Assume first that 1) holds, and that there is an exact sequence of $R\C$-modules
\[
0\rightarrow L \rightarrow N\rightarrow M \rightarrow 0.
\]
Since the functor $\res_S$ is exact and the functor $\ind_S$ right exact, we get a commutative diagram with exact rows
\[
\xymatrix{
& \ind_S\res_S L \ar[d]_{\mu_L} \ar[r] & \ind_S\res_S N \ar[d]_{\mu_N} \ar[r] & \ind_S\res_S M \ar[d]_{\mu_M} \ar[r] & 0 \ar[d]\\
0\ar[r]& L \ar[r] &  N  \ar[r]  &  M \ar[r] & 0 }
\]
where $\mu_M$ is an isomorphism and $\mu_N$ is an epimorphism. An easy diagram chase shows us that $\mu_L$ is an epimorphism, so 2) holds.

The implication 2) $\Rightarrow$ 3) is trivial. Assume next that 3) holds. Now the morphism $\bigoplus_{s\in S} K(s)[\Mor_{\C}(s,-)]\rightarrow K$ is an epimorphism, so both $\bigoplus_{s\in S} K(s)[\Mor_{\C}(s,-)]$ and $K$ have the same image in $\bigoplus_{s\in S} M(s)[\Mor_{\C}(s,-)]$. The required exactness then follows immediately.

Trivially 4) implies 5). Finally, let us assume that 5) holds. By Proposition \ref{induktio_apu}, we get a commutative diagram with exact rows
\[
\xymatrix{
\bigoplus_{s\in S} B_s[\Mor_{\C}(s,-)] \ar[d]_{\cong} \ar[r] & \bigoplus_{s\in S} A_s[\Mor_{\C}(s,-)]  \ar[d]_{\cong} \ar[r] & \ind_S\res_S M \ar[d]_{\mu_M} \ar[r] & 0 \ar[d]\\
\bigoplus_{s\in S} B_s[\Mor_{\C}(s,-)] \ar[r] &  \bigoplus_{s\in S} A_s[\Mor_{\C}(s,-)]  \ar[r]  &  M \ar[r] & 0 }
\]
from which we can see that $\mu_M$ is an isomorphism by the five lemma.
\end{proof}

\begin{remark}\label{djament_huom}
Proposition~\ref{s-esitetyt} is due to Djament~\cite{Djament}*{p.~11, Prop.~2.14}. The reader should be cautious, since we use the term `support' in a different meaning as in \cite{Djament}.
\end{remark}

The following proposition is a special case of~\cite{Gabriel}*{p.~83, Prop.}. For the sake of clarity, we present a proof using our notation.

\begin{proposition}\label{roiter_theorem}
An $R\C$-module $M$ is finitely presented if and only if there exists a finite full subcategory $S\subseteq \C$ such that
\begin{itemize}
\item[1)] $M(s)$ is finitely presented for all $s\in S$;
\item[2)] $M$ is $S$-presented.
\end{itemize}
\end{proposition}

\begin{proof}
Assume first that $M$ is finitely presented, so that there exists an exact sequence
\[
\bigoplus_{j\in J} R[\Mor_{\C}(b_j,-)] \rightarrow \bigoplus_{i\in I}R[\Mor_{\C}(a_i,-)] \rightarrow M\rightarrow 0,
\]
where $I$ and $J$ are finite sets, and $a_i, b_j\in \C$ for all $i\in I$ and $j\in J$. Evaluating this at point $c\in \C$ gives us an exact sequence
\[
R^{m_c} \rightarrow R^{n_c} \rightarrow M(c)\rightarrow 0
\]
for some $m_c, n_c\in \mathbb N$, so that 1) holds. For 2), by setting
\[
S:=\{a_i\mid i\in I\} \cup \{b_j\mid j\in J\}
\]
we immediately see that $M$ is $S$-presented by Proposition~\ref{s-esitetyt} 5).

Assume next that there exists a finite full subcategory $S\subseteq \C$ such that 1) and 2) hold. Now $M$ is $S$-generated, so the natural morphism $\bigoplus_{s\in S}M(s)[\Mor_\C(s,-)]\rightarrow M$ is an epimorphism. Since $M(s)$ is finitely generated for all $s\in S$, there exists an epimorphism $R^{n_s}\rightarrow M(s)$ for all $s\in S$, where $n_s\in \mathbb N$. Combining these epimorphisms, we get an epimorphism
\[
\bigoplus_{t\in S}R^{n_t}[\Mor_\C(t,-)]\rightarrow \bigoplus_{t\in S}M(t)[\Mor_\C(t,-)] \rightarrow M
\] 
and an exact sequence
\[
0\rightarrow N \rightarrow \bigoplus_{t\in S}R^{n_t}[\Mor_\C(t,-)]\rightarrow M \rightarrow 0.
\]
Because $M$ is $S$-presented, $N$ must be $S$-generated by Proposition~\ref{s-esitetyt} 2), so there exists an epimorphism $\bigoplus_{t\in S}N(t)[\Mor_S(t,-)]\rightarrow N$. On the other hand, $M(s)$ is finitely presented, so $N(s)$ is finitely generated for all $s\in S$. Thus there exists an epimorphism $R^{m_s}\rightarrow N(s)$ for all $s\in S$, where $m_s\in \mathbb N$. Hence we get an exact sequence
\[
\bigoplus_{t\in S}R^{m_t}[\Mor_\C(t,-)]\rightarrow \bigoplus_{t\in S}R^{n_t}[\Mor_\C(t,-)]\rightarrow M\rightarrow 0.
\]
\end{proof}

From the proof of Proposition~\ref{roiter_theorem} we immediately get the following corollary:

\begin{corollary}\label{fin_gen}
An $R\C$-module $M$ is finitely generated if and only if there exists a finite full subcategory $S\subseteq \C$ such that
\begin{itemize}
\item[1)] $M(s)$ is finitely generated for all $s\in S$;
\item[2)] $M$ is $S$-generated.
\end{itemize}
\end{corollary}


\subsection{Births and deaths relative to \texorpdfstring{$S$}{S}}

From now on, we will assume that $\C$ is a poset. 

Let $M$ be an $R\C$-module, $S\subseteq \C$ a subset, and $c\in \C$. Write $S':=S\backslash \{c\}$. We note that
\[
\colim_{d < c, \ d\in S} M(d) = \colim_{d\leq c, \ d\in S'} M(d) = (\ind_{S'}\res_{S'} M)(c).
\]
Since $\res_{S'}$ is exact and $\ind_{S'}$ is right exact, we then see that the functor
\[
\RCmod \rightarrow \Rmod, \ M \mapsto \colim_{d<c, \ d\in S} M(d).
\]
is also right exact.

\begin{definition}\label{syntyma_kuolema}
Let $\C$ be a poset, $M$ an $R\C$-module, $S\subseteq \C$ a subset and $c\in \C$. Let 
\[\lambda_{M,c}\colon  \colim_{d<c, \ d\in S} M(d) \rightarrow M(c)\]
be the natural homomorphism. We define the \emph{set of births relative to $S$} by
\[
B_S(M):=\{c\in \C\mid \lambda_{M,c} \ \text{is a non-epimorphism}\}
\]
and the \emph{set of deaths relative to $S$} by
\[
D_S(M):=\{c\in \C\mid \lambda_{M,c} \ \text{is a non-monomorphism}\}.
\] 
\end{definition}

\begin{remark}\label{syntyma_inkluusio}
Note that $\lambda_{M,c}$ is an epimorphism if and only if the natural homomorphism $\bigoplus_{d<c, \ d\in S} M(d)\rightarrow M(c)$ is an epimorphism. This implies that if $T\subseteq S\subseteq \C$, then $B_S(M)\subseteq B_T(M)$.
\end{remark}

\begin{example}\label{interval_ex}
Let $\C$ be a poset. Let $I$ be an \emph{interval} of $\C$ i.e.~a non-empty subset of $\C$ satisfying the condition that if $a,b\in I$, $c\in \C$ and $a\leq c\leq b$, then $c\in I$. Let $R_I$ be the $R\C$-module defined on objects by
\[
R_I(c)=\left\{ \begin{aligned} 
&R, \ \text{when} \ c\in I;\\
&0, \ \text{otherwise,}
\end{aligned}\right.
\]
and with identity morphisms inside the interval. Then the sets of births $B_{\C}(R_I)$ and $B_{I}(R_I)$ both consist of the minimal elements of $I$. To find the deaths, we note that $R_I$ is ${\uparrow} I$-presented, so deaths must either be inside $I$ or above it (see Remark \ref{kakkosehto}).

First, let $c\in S_1:=({\uparrow I})\backslash I$. Now $R_I(c)=0$. Since $S_1\subseteq {\uparrow}\supp(R_I)$, we see that $\colim_{d<c, \ d\in I}R_I(d) \neq 0$. Thus $c\in D_I(R_I)$, so $S_1\subseteq D_I(R_I)$. Furthermore, it is clear that $\colim_{d<c}R_I(d)\neq 0$ if and only if $c$ is minimal in $S_1$. This implies that exactly the minimal elements of $S_1$ are in $D_\C(R_I)$.

Secondly, let $c\in I$. It is straightforward to see that $c\in D_I(R_I)$ if and only if the set $(I\cap{\downarrow}c) \backslash \{c\}$ is not connected as a poset. This applies also  
to $D_\C(R_I)$. Set
\[
S_2 := \{c\in I \mid (I\cap{\downarrow}c) \backslash \{c\} \text{ is not connected }\}.
\]
We conclude that $D_I(R_I)=S_1\cup S_2$, while $D_\C(R_I)$ is the union of the set of the minimal elements of $S_1$, and the set $S_2$.
\end{example}

\begin{proposition}\label{FSP_apu}
Let $\C$ be a poset, $M$ an $R\C$-module, and $S\subseteq \C$ a subset. Then
\begin{itemize}
\item[1)] $M$ is $S$-generated if and only if $B_S(M)\subseteq S$;
\item[2)] $M$ is $S$-presented if and only if $B_S(M)\cup D_S(M) \subseteq S$.
\end{itemize}
\end{proposition}

\begin{proof}
Both 1) and 2) are proved similarly. We only prove 2) here. Directly from the definitions,
\begin{align*}
&M \text{ is }S\text{-presented} \\ &\Leftrightarrow \ \mu_{M,c}\colon  \colim_{d\leq c, \ d\in S} M(d)\rightarrow M(c) \text{ is an isomorphism for all }c\in \C\\
&\Leftrightarrow \ \lambda_{M,c} \colon \colim_{d<c, \ d\in S} M(d)\rightarrow M(c) \text{ is an isomorphism for all }c\in \C\backslash S\\ 
&\Leftrightarrow \ B_S(M)\cup D_S(M)\subseteq S.
\end{align*}
\end{proof}

Let $S\subseteq \C$. We may think of $B_\C(M)$ as the set of `real' births of $M$. The following proposition shows that for $S$-generated modules we may focus only on births relative to $S$.

\begin{proposition}\label{syntyma_vertailu}
Let $M$ be an $S$-generated $R\C$-module. Then $B_\C(M)=B_S(M)$.
\end{proposition}

\begin{proof}
By Remark~\ref{syntyma_inkluusio} it is enough to show that $B_S(M)\subseteq B_{\C}(M)$. Let $c\in \C\backslash B_\C(M)$, so that the natural homomorphism $\bigoplus_{d<c}M(d)\rightarrow M(c)$ is an epimorphism. Since $M$ is $S$-generated, there is an epimorphism $\bigoplus_{d'\leq d, \ d'\in S}M(d')\rightarrow M(d)$ for all $d<c$. We may combine these epimorphisms to get an epimorphism $\bigoplus_{d<c, \ d\in S}M(d)\rightarrow M(c)$, implying that $c\in \C\backslash B_S(M)$.
\end{proof}


\subsection{\texorpdfstring{$S$}{S}-splitting}

\begin{definition}\label{s-lohkaisu}
Let $S\subseteq \C$ a subset and $c\in \C$. If $M$ is an $R\C$-module, denote by $S_{S,c}M$ the $R$-module defined by the exact sequence 
\[
\colim_{d<c, \ d\in S}M(d) \stackrel{\lambda_{M,c}}{\rightarrow} M(c) \stackrel{\pi_{M,c}}{\rightarrow} S_{S,c}M \rightarrow 0,
\]
where $\pi_{M,c}$ is the canonical epimorphism. This gives rise to a functor $S_{S,c}$, the \emph{$S$-splitting functor at $c$}. More explicitly,
\[
S_{S,c}M = M(c) / \Img(\lambda_{M,c}).
\]
\end{definition}

\begin{remark}\label{lohkaisu_eksakti}
The $S$-splitting functor at $c$ could equivalently be defined as the composition of the splitting functor $S_c$~(\cite{Luck}*{p. 156}) and the restriction functor $\res_{S\cup \{c\}}$ by setting $S_{S,c}=S_c\circ \res_{S\cup\{c\}}$. Since both $S_c$ and $\res_{S\cup \{c\}}$ are left adjoints, we see that $S_{S,c}$ is a left adjoint, and thus additive.
\end{remark}

The concept of a splitting functor is due to Lück (\cite{Luck}). The basic example for us is the following:

\begin{example}\label{split_esim}
Let $k$ be a field, $S\subseteq \mathbb Z^n$ a subset, and $M$ a $k(\mathbb N^n{\smallint} \mathbb Z^n)$-module. We identify $M$ with the corresponding $\mathbb Z^n$-graded $k[X_1,\ldots,X_n]$-module. Denote by $m:=\langle X_1,\ldots,X_n\rangle$ the maximal homogeneous ideal of $k[X_1,\ldots,X_n]$. If $N$ is the homogeneous submodule of $M$ generated by the union of $M_s$, where $s\in S$, we notice that
\[
(M/mN)_c =  M_c / (mN)_c = M(c) / \Img(\lambda_{M,c}) = S_{S,c} M
\]
for all $c\in\mathbb Z^n$. In particular, this yields an isomorphism of $k$-vector spaces, 
\[
M/mN \cong \bigoplus_{c\in \mathbb Z^n} S_{S,c}M.
\]
\end{example}

\begin{remark}\label{lohkaisu_syntyma}
Let $M$ be an $R\C$-module and $S\subseteq \C$ a subset. Note that for all $c\in \C$, we have $c\in B_S(M)$ if and only if $S_{S,c} M \neq 0$.
\end{remark}

\begin{remark}\label{colim_apu}
Let $A$ be an $R$-module, $S\subseteq \C$ a subset, $s\in S$, and $c\in \C$. Let $S':=S\backslash\{c\}$. If $s\neq c$, we see that
\[
\colim_{d < c, \ d\in S} A[\Mor_\C(s,d)] =\colim_{d\leq c, \ d\in S'} A[\Mor_\C(s,d)] \cong A[\Mor_\C(s,c)]
\]
by Proposition \ref{induktio_apu}. If $s=c$, then obviously $\colim_{d<c, \ d\in S} A[\Mor_\C(s,d)]=0$. In particular
\[
S_{S,c}(A[\Mor_\C(s,-)])=\left\{\begin{aligned} &A, \text{ when }s=c\\
&0, \ \text{otherwise}.
\end{aligned} \right.
\]
\end{remark}

Next, we prove a version of Nakayama's lemma (cf.~\cite{Tchernev}*{p.~12, Lemma~6.2}).

\begin{lemma}\label{nakayama}
Let $M$ be an $R\C$-module and $S\subseteq \C$ a subset. If $\supp(M)\cap S$ has a minimal element $c$, then $S_{S,c}M\neq 0$.
\end{lemma}

\begin{proof}
Assume that $c\in\supp(M)\cap S$ is minimal. Then $M(d)=0$ for all $d\in S$ with $d<c$. In particular, $\colim_{d<c, \ d\in S} M(d) = 0$. Thus $S_{S,c} M\neq 0$.
\end{proof}

Recall that a poset $P$ is called \emph{Artinian}, if there are no infinite strictly descending chains of elements of $P$, or equivalently, if every non-empty subset $S\subseteq P$ has a minimal element.

\begin{proposition}\label{tchernev_63}
Let $f\colon L\rightarrow M$ be a morphism of $R\C$-modules, where $M$ is $S$-generated with an Artinian $S\subseteq \C$. If $S_{S,c}f\colon  S_{S,c}L\rightarrow S_{S,c}M$ is an epimorphism for all $c\in B_S(M)$, then $f$ is an epimorphism.
\end{proposition}

\begin{proof}
We first note that $\Coker f$ is $S$-generated, since $M$ is $S$-generated. Suppose that $f$ is not an epimorphism. Then $\Coker f\neq 0$, so there exists $s\in S$ such that $(\Coker f)(s)\neq 0$. Hence $\supp(\Coker f)\cap S$ has a minimal element $c$ by the Artinian property. Now $S_{S,c}(\Coker f)\neq 0$ by Lemma~\ref{nakayama}, which implies that $c\in B_S(\Coker f)\subseteq B_S(M)$. Since $S_{S,c}$ is right exact, we get $\Coker S_{S,c}f \neq 0$, so $S_{S,c}f$ is not an epimorphism. 
\end{proof}

\begin{lemma}\label{sf-lemma}
Let $S\subseteq \C$ be a subset and
\[
0\rightarrow L \stackrel{j}{\rightarrow}N\stackrel{f}{\rightarrow}M\rightarrow 0
\]
an exact sequence of $R\C$-modules. The following are equivalent for all $c\in \C$:
\begin{itemize}
\item[1)] $(\Ker f)(c) \subseteq  \Img \lambda_{N,c}$;
\item[2)] $S_{S,c}(j)=0$;
\item[3)] $S_{S,c}(f)$ is a monomorphism;
\item[4)] $S_{S,c}(f)$ is an isomorphism.
\end{itemize}
\end{lemma}

\begin{proof}
The equivalence of $1)$ and $3)$ immediately follows from the fact that
\[
\Ker S_{S,c}(f) = ((\Ker f)(c)+\Img \lambda_{N,c})/\Img \lambda_{N,c}.
\]
Since $S_{S,c}$ is right exact, we have $\Ker S_{S,c}(f) = \Img S_{S,c}(j)$. Therefore $2)$ implies $3)$. The equivalence of $3)$ and $4)$ holds, because $S_{S,c}$ preserves epimorphisms.
\end{proof}

We recall that an epimorphism of $R\C$-modules $f\colon N\rightarrow M$ is called \emph{minimal}, if for all morphisms $g\colon L\rightarrow N$, $fg$ is an epimorphism if and only if $g$ is an epimorphism.
It is known that an epimorphism $f$ is minimal if and only if for all \hyphenation{sub-modules} submodules $N'\subseteq N$
\[
N'+\Ker f=N \ \Rightarrow \ N'=N.
\] 
A minimal epimorphism $f\colon N\rightarrow M$, where $N$ is projective, is called a \emph{projective cover} of $M$ (see e.g.~\cite{Assem}*{p.~28}).

\begin{remark}\label{proj_huomautus}
Let $A$ be an $R$-module and $c\in \C$. Then $A$ may be thought of as an $R\{c\}$-module, and we note that $A[\Mor_\C(c,-)] \cong \ind_{\{c\}}A$. In particular, the functor $A \mapsto A[\Mor_\C(c,-)]$ preserves projectives, since it is the left adjoint of the exact functor $\res_{\{c\}}$.
\end{remark}

\begin{proposition}\label{turhuuslemma}
Let $f\colon N\rightarrow M$ be an epimorphism of $S$-generated $R\C$-modules, where $S\subseteq \C$ is Artinian. If $(\Ker f)(c)\subseteq \Img \lambda_{N,c}$ for all $c\in \C$, then $f$ is minimal.
The converse implication holds if $S_{S,c}M$ is projective for all $c\in S$. 
\end{proposition}

\begin{proof}
Let $(\Ker f)(c)\subseteq \Img \lambda_{N,c}$ for all $c\in\C$. Suppose that $N'+\Ker f=N$ for some submodule $N'\subseteq N$. We note that for all $c\in \C$, $(N')(c)+\Img\lambda_{N,c} = N(c)$. This implies that $S_{S,c}N' = S_{S,c}N$ for all $c\in \C$. Since $S$ is Artinian, we may use Proposition~\ref{tchernev_63} to conclude that $N'=N$, so $f$ is minimal.

Next, let $f$ be minimal, and let $S_{S,c}M$ be projective for all $c\in S$. Thus we can find sections $S_{S,c}M\rightarrow M(c)$ for all $c\in S$. These induce a morphism
\[
h\colon  \bigoplus_{s\in S} S_{S,s}M[\Mor_\C(s,-)] \rightarrow \bigoplus_{s\in S} M(s)[\Mor_\C(s,-)] \rightarrow M.
\]
Remark~\ref{colim_apu} now implies that $S_{S,c}h=\id_{S_{S,c}M}$ for all $c\in S$, so $h$ is an epimorphism by Proposition~\ref{tchernev_63}. 

Since $S_{S,c}M$ is projective for all $c\in S$, we see that $\bigoplus_{s\in S} S_{S,s}M[\Mor_\C(s,-)]$ is also projective by Remark~\ref{proj_huomautus} (as a sum of projectives). Thus the morphism $h$ factors through $f$, and we get a diagram
\[
\xymatrix{
\bigoplus_{s\in S}S_{S,s}M[\Mor_c(s,-)] \ar[r]^-{g} \ar[dr]_-{h} & N\ar[d]^{f}\\
& M }
\]
that commutes. Now $f$ is minimal, so $g$ is an epimorphism. Applying functor $S_{S,c}$, where $c\in S$, on the diagram, we see that $S_{S,c}f\circ S_{S,c}g=\id$, which implies that $S_{S,c}g$ is a monomorphism, and therefore an isomorphism. Hence $S_{S,c}f$ is an isomorphism for all $c\in S$. This is equivalent to $(\Ker f)(c) \subseteq \Img \lambda_{N,c}$ for all $c\in \C$ by Lemma~\ref{sf-lemma}.
\end{proof}

\begin{remark}\label{proj_verho}
Let $M$ be an $S$-generated $R\C$-module, where $S$ is Artinian. If $S_{S,c}M$ is projective for all $c\in S$, the morphism $h\colon  \bigoplus_{s\in S}S_{S,s}M[\Mor_c(s,-)] \rightarrow M$ induced by sections $S_{S,s}M\rightarrow M(s)$ is a projective cover of $M$.

Indeed, as noted earlier, $h$ is an epimorphism with $S_{S,c}h=\id$ for all $c\in S$. Then Lemma~\ref{sf-lemma} implies that $(\Ker f)(c) \subseteq \Img \lambda_{N,c}$ for all $c\in\C$, and the rest follows from Proposition~\ref{turhuuslemma}.
\end{remark}


\subsection{Minimality of births and deaths}

We will now show how the sets of births and deaths relative to a subset $S\subseteq \C$ are in a sense minimal if the module is $S$-generated or $S$-presented.

\begin{proposition}\label{syntyma_minimi}
Let $M$ be an $S$-generated $R\C$-module, where $S\subseteq \C$ is Artinian. Then $M$ is $B_S(M)$-generated. Furthermore, $B_S(M)$ is the minimum element of the set $\{T\subseteq S\mid M \text{ is }T\text{-generated}\}$.
\end{proposition}

\begin{proof}
Let $\rho$ be the natural morphism $\rho\colon  \bigoplus_{s\in B_S(M)}M(s)[\Mor_\C(s,-)] \rightarrow M$. Remark~\ref{colim_apu} shows us that applying the $S$-splitting functor at $c\in S$ yields the canonical epimorphism $S_{S,c}\rho=\pi\colon  M(c)\rightarrow S_{S,c}M$. Thus $\rho$ is an epimorphism by Proposition~\ref{tchernev_63}.

To show the claimed minimality: If $M$ is $T$-generated for some $T\subseteq S$, we have $B_S(M)\subseteq B_T(M) \subseteq T$ by Proposition~\ref{FSP_apu} and Remark~\ref{syntyma_inkluusio}.
\end{proof}

Next, we introduce a technical lemma.

\begin{lemma}\label{kaarme}
Assume that we have a commutative diagram of $R$-modules with exact rows
\[
\xymatrix{
& L\ar[d]_{f} \ar[r] & N \ar[d]_{g} \ar[r] & M \ar[d]^{h}\ar[r]& 0 \\
0 \ar[r] & L' \ar[r] & N' \ar[r] & M'& }
\]
where $g$ is a monomorphism. If $f$ is an epimorphism, then $h$ is a monomorphism. The converse holds if either the natural morphism $\Coker g \rightarrow \Coker h$ is a monomorphism or $g$ is an epimorphism. 
\end{lemma}

\begin{proof}
The snake lemma gives us an exact sequence
\[
\Ker f \rightarrow  \Ker g \rightarrow  \Ker h \rightarrow  \Coker f \rightarrow  \Coker g \rightarrow\Coker h,
\]
where $\Ker g=0$. If $\Coker f=0$, we get $\Ker h=0$. If $\Coker g=0$, we have $\Ker h\cong\Coker f$, and we are done. If $\Coker g \rightarrow \Coker h$ is a monomorphism, we see that $\Coker f$ maps to $0$, so $\Ker h \rightarrow \Coker f$ is an epimorphism. Since $\Ker g=0$, the morphism $\Ker h \rightarrow \Coker f$ is also a monomorphism.
\end{proof}

\begin{lemma}\label{esitys_minimi}
Let $M$ be an $S$-presented $R\C$-module, where $S\subseteq \C$ is Artinian. Assume that we have an exact sequence of $R\C$-modules
\[
0\rightarrow L \rightarrow N \stackrel{f}{\rightarrow} M\rightarrow 0,
\]
where $N$ is $S$-generated and $D_S(N)=\emptyset$. Then $D_S(M)\subseteq B_S(L)$. Furthermore, if $N$ is $B_S(M)$-generated, we have $B_S(L)\subseteq B_S(M)\cup D_S(M)$.
\end{lemma}

\begin{proof}
Let $c\in \C$. Applying $\colim_{d< c, \ d\in S}$ to the exact sequence above, we get a diagram with exact rows
\[
\xymatrix@C=1.5em{
& \colim_{d<c, \ d\in S}L(d)\ar[d]_{\lambda_{L,c}} \ar[r] & \colim_{d<c, \ d\in S}N(d) \ar[d]_{\lambda_{N,c}} \ar[r] & \colim_{d<c, \ d\in S}M(d)\ar[r] \ar[d]_{\lambda_{M,c}}& 0 \\
0 \ar[r] & L(c)\ar[r] & N(c) \ar[r]^-{f_c} & M(c) \ar[r]& 0 }
\]
that commutes. Here $\lambda_{N,c}$ is a monomorphism, because $D_S(N)=\emptyset$. To show that $D_S(M)\subseteq B_S(L)$, suppose that $c\notin B_S(L)$. In this case $\lambda_{L,c}$ is an epimorphism, so $\lambda_{M,c}$ is a monomorphism by Lemma~\ref{kaarme}. Thus $c\notin D_S(M)$.

Assume that $N$ is $B_S(M)$-generated. Since $f$ is an epimorphism, Proposition~\ref{syntyma_minimi} implies that $B_S(M)=B_S(N)$. Suppose that $c\notin B_S(M)\cup D_S(M)$. Now $\lambda_{N,c}$ is an epimorphism, since $c\notin B_S(N)=B_S(M)$. Moreover, $\lambda_{M,c}$ is a monomorphism because $c\notin D_S(M)$. It follows from Lemma~\ref{kaarme} that $\lambda_{L,c}$ \text{ is an epimorphism}, so $c\notin B_S(L)$. Thus $B_S(L)\subseteq B_S(M)\cup D_S(M)$.
\end{proof}

\begin{proposition}\label{esitys_minimi_osa2}
Let $M$ be an $S$-presented $R\C$-module, where $S\subseteq \C$ is Artinian. Then $M$ is $B_S(M)\cup D_S(M)$-presented. Furthermore, $B_S(M)\cup D_S(M)$ is the minimum element of the set $\{T\subseteq S\mid M \text{ is }T\text{-presented}\}$.
\end{proposition}

\begin{proof}
Let us examine an exact sequence
\[
0\rightarrow L \rightarrow N \rightarrow M\rightarrow 0,
\]
with $N$ of the form $N=\bigoplus_{s\in B_S(M)} A_s[\Mor_\C(s,-)]$, where $A_s$ is an $R$-module for all $s\in B_S(M)$. Note that such $N$ always exists by Proposition~\ref{syntyma_minimi}. Since $M$ is $S$-presented, Proposition~\ref{s-esitetyt} 2) implies that $L$ is $S$-generated. Using Proposition~\ref{s-esitetyt} 5), we notice that if $L$ is $T$-generated for some $T\subseteq S$, then $M$ is $(B_S(M)\cup T)$-presented. Now $L$ is $B_S(L)$-generated by Proposition~\ref{syntyma_minimi}, so we deduce that $M$ is $(B_S(M)\cup B_S(L))$-presented. We can now use Lemma~\ref{esitys_minimi} to see that then $M$ is$(B_S(M)\cup D_S(M))$-presented.

Suppose next that $M$ is also $T$-presented for some $T\subseteq S$. As in the proof of Lemma~\ref{esitys_minimi}, we note that $B_S(M)=B_S(N)$. The minimality of $B_S(M)$ in Proposition~\ref{syntyma_minimi} implies that $B_S(N)=B_S(M)\subseteq T$, so $N$ is $T$-generated by Proposition~\ref{syntyma_minimi}. Thus $L$ is $T$-generated by Proposition~\ref{s-esitetyt} 2). Therefore we must have 
\[
B_S(M)\subseteq B_T(M)\subseteq T \ \text{and} \ B_S(L)\subseteq B_T(L)\subseteq T,
\]
by Proposition~\ref{FSP_apu} 1) and Remark~\ref{syntyma_inkluusio}. We use Lemma~\ref{esitys_minimi} to conclude that
\[
B_S(M) \cup D_S(M) = B_S(L)\cup B_S(M) \subseteq T.
\]
\end{proof}

\begin{remark}\label{verho_lemma}
Assume that $M$ is an $S$-presented $R\C$-module, where $S\subseteq \C$ is Artinian. Let $f\colon N\rightarrow M$ be a projective cover. Then $S_{S,c}f$ is an isomorphism for all $c\in \C$ if and only if $S_{S,c}M$ is projective for all $c\in \C$.

To see this, first suppose that $S_{S,c}f$ is an isomorphism for all $c\in \C$. Since $S_{S,c}$ preserves projectives for all $c\in\C$, we see that $S_{S,c}N$ is projective, and thus $S_{S,c}M$ is projective.

Conversely, suppose that $S_{S,c}M$ is projective for all $c\in \C$. We may now apply Proposition~\ref{turhuuslemma} and Lemma~\ref{sf-lemma} to get isomorphisms $S_{S,c}f\colon  S_{S,c}N\rightarrow S_{S,c}M$ for all $c\in \C$.
\end{remark}

\begin{remark}\label{rank_huom} 
In~\cite{Carlsson}, Carlsson and Zomorodian define multiset-valued invariants $\xi_0$ and $\xi_1$ for a finitely generated $\mathbb Z^n$-graded $k[X_1,\ldots,X_n]$-module $M$, where $k$ is a field. 
The multisets $\xi_0(M)$ and $\xi_1(M)$ indicate the degrees in $\mathbb Z^n$ where the elements of $M$ are born and where they die, respectively. In more algebraic terms, $\xi_0(M)$ and $\xi_1(M)$ consist
of the degrees of minimal generators and minimal relations of M equipped with the multiplicities they occur. Consider an exact sequence
\[
0\rightarrow L \rightarrow N \stackrel{f}{\rightarrow} M \rightarrow 0,
\]
where $N$ is a free module and $f$ a minimal homomorphism. Since $M$ is $S$-presented for some finite $S\subseteq \mathbb Z^n$, it is easy to see that $\xi_0(M)$ is a multiset where the underlying set is $B_S(M)$ and the multiplicity of $c\in B_S(M)$ is the dimension of $M(c)$. Note that the choice of $S$ does not matter here, since $B_S(M)=B_\C(M)$ by Proposition \ref{syntyma_vertailu}. We note that $L$ is $S$-generated by Proposition \ref{s-esitetyt} 3), so we may apply a similar argument to conclude that $\xi_1(M)$ is a multiset with $B_S(L)$ as the underlying set and the dimension of $L(c)$ as the multiplicity of $c\in B_S(L)$. The next theorem will show that $D_S(M)$ is the underlying set of $\xi_1(M)$.
\end{remark}

\begin{theorem}\label{verho}
Let $M$ be an $S$-presented $R\C$-module, where $S\subseteq \C$ is Artinian. Assume that $S_{S,c}M$ is projective for all $c\in B_S(M)$. If
\[
0\rightarrow L \rightarrow N\stackrel{f}{\rightarrow} M\rightarrow 0
\]
is an exact sequence where $f$ is a projective cover, then $B_S(L)= D_S(M)$.
\end{theorem}

\begin{proof} By Lemma~\ref{esitys_minimi}, it is enough to show that $B_S(L) \subseteq D_S(M)$. Let $c\in B_S(M)$.
Suppose that we would have $c\notin D_S(M)$. Then $\lambda_{M,c}$ is a monomorphism. Since $f$ is minimal, by Proposition~\ref{turhuuslemma} and Lemma~\ref{sf-lemma} there exists a natural isomorphism
\[
S_{S,c}f\colon  S_{S,c}N=\Coker \lambda_{N,c} \rightarrow \Coker \lambda_{M,c}=S_{S,c}M.
\]
It now follows from Lemma~\ref{kaarme} that $\lambda_{L,c}$ is an epimorphism, which is equivalent to $c\notin B_S(L)$. But this is impossible.
\end{proof}


\section{Presentations with finite support}

In this section we will prove our main result, Theorem~\ref{fin_pres}, which gives a characterization for finitely presented modules. We will assume that $\C$ is a poset and $R$~a~commutative ring. 


\subsection{\texorpdfstring{$S$}{S}-determined \texorpdfstring{$R\C$}{RC}-modules}

Let $M$ be an $R\C$-module. If $S\subseteq \C$ is a finite set such that $M$ is $S$-presented, we say that $S$ is a \emph{finite support of a presentation (FSP)} of $M$. In what follows, we are trying to find a  condition equivalent for $M$ having an FSP.

\begin{definition}\label{s-determined}
An $R\C$-module $M$ is \emph{$S$-determined} if there exists a subset $S\subseteq \C$ such that $\supp(M)\subseteq {\uparrow} S$, and for every $c\leq d$ in $\C$
\[
S\cap {\downarrow} c =  S \cap {\downarrow} d \ \Rightarrow \ M(c\leq d) \text{ is an isomorphism.}
\]
\end{definition}

\begin{remark}\label{kakkosehto}
Let $M$ be an $R\C$-module and $S\subseteq \C$ a subset. Denote $T:={\uparrow} S$. Then the condition $\supp(M)\subseteq T$ of Definition~\ref{s-determined} is equivalent to the following conditions:
\begin{itemize}
\item[1)] $M$ is $T$-generated;
\item[2)] $M$ is $T$-presented;
\item[3)] If $S\cap{\downarrow} c = \emptyset$, then $M(c)=0$.
\end{itemize}
To show this, we first note that $1)$ implies $3)$, because ${\uparrow}S=T$. Taking the contraposition of 3), we get $\supp(M)\subseteq T$. Next, note that below every $c\in D_T(M)$ there must be some $d~\in~\supp(M)$ such that $d<c$. Thus $D_T(M)\subseteq {\uparrow} \supp(M)$. Obviously also $B_T(M)\subseteq \supp(M)$. We now observe that if $\supp(M)\subseteq T$, we get
\[
B_T(M)\cup D_T(M)\subseteq {\uparrow} \supp(M) \subseteq {\uparrow} T = T,
\]
This means that $\supp(M)\subseteq T$ implies 2) by Proposition~\ref{FSP_apu}. Finally, 3) trivially follows from 2).
\end{remark}

\begin{proposition}\label{FSP_ehdot}
Let $M$ be an $S$-presented $R\C$-module, where $S\subseteq \C$. Then $M$ is $S$-determined.
\end{proposition}

\begin{proof}
Trivially $\supp(M)\subseteq {\uparrow} S$. If $c\leq d$ in $\C$, we have a commutative diagram
\[
\xymatrix{
\colim_{e\leq c, \ e\in S}M(e)\ar[d]_{\cong} \ar[r] & \colim_{e\leq d, \ e\in S}M(e)  \ar[d]^{\cong} \\
M(c) \ar[r]_{M(c\leq d)} &  M(d)}
\]
with the vertical isomorphisms being components of the canonical isomorphism of Proposition~\ref{s-esitetyt}, 1). This immediately shows us that $M$ is $S$-determined. 
\end{proof}


\subsection{Minimal upper bounds}

Let $S\subseteq \C$. We would like to find conditions under which $S$-determined implies $S$-presented. In general this is false (see Example~\ref{mub_esimerkki}), so we first need to apply some technical limitations on the poset $\C$ to guarantee that it is "small" enough.

\begin{notation}
Let $S\subseteq \C$ be a subset. An element $c\in \C$ is an \emph{upper bound of $S$}, if $c\geq s$ for all $s\in S$. We denote the set of minimal upper bounds of $S$ by $\mub_\C(S)$.
\end{notation}

\begin{definition}
The poset $\C$ is \emph{weakly bounded from above} if every finite $S\subseteq \C$ has a finite number of minimal upper bounds in $\C$.
\end{definition}

\begin{definition}
The poset $\C$ is \emph{mub-complete} if given a finite non-empty subset $S\subseteq\C$ and an upper bound $c$ of $S$, there exists a minimal upper bound $s$ of $S$ such that $s\leq c$. 
\end{definition}

\begin{remark}\label{huom8}
A poset that is weakly bounded from above and mub-complete is called a poset with property $\mathcal M$ in~\cite{Lawson}. In contrast to~\cite{Lawson}, we do not require the empty set to have minimal upper bounds for a poset to be mub-complete.
\end{remark}

\begin{example}
Let $L$ be a poset where every finite subset $K\subseteq L$ has an infimum and a supremum (i.e.~$L$ is a lattice). Then $L$ is weakly bounded from above and mub-complete.
\end{example}

A `good' monoid $G$ in~\cite{Corbet} is a cancellative monoid that is weakly bounded from above as a poset (with the natural order). If $G$ is also commutative, we get the following description of mub-completeness.

\begin{proposition}\label{mub_mcd}
Let $G$ be a commutative cancellative monoid that is weakly \hyphenation{boun-ded} bounded from above as a poset (with the natural order). Then $G$ is mub-complete if and only if there exists a maximal common divisor for each $g,h\in G$.
\end{proposition}

\begin{proof}
Assume first that $G$ is mub-complete. Let $g,h\in G$. Since $gh$ is an upper bound of $g$ and $h$, there exists a minimal upper bound $j\in G$ of $g$ and $h$ such that $lj=gh$ for some $l\in G$.  We claim that $l$ is a maximal common divisor of $g$ and $h$. We may write $j=ag=bh$, where $a,b\in G$. Now
\[
gh=lj=lag=lbh,
\] 
so that $g=lb$ and $h=la$ by cancellativity. Thus $l$ is a common divisor of $g$ and $h$. For the maximality, let $k\in G$ be another common divisor of $g$ and $h$ such that $l$ divides $k$. We may then write $k=k'l$, where $k'\in G$. Furthermore, we have $g=ck$ and $h=dk$ for some $c,d\in G$. Combining these equations, we get
\[
lj=gh=ck'lh=gdk'l.
\]
Cancelling $l$, we see that $j=k'ch=k'dg$. Furthermore, cancelling $k'$ yields $ch=dg$, another upper bound for $g$ and $h$. Since $j$ is a minimal upper bound of $g$ and $h$, we must have $k'=1$, proving the maximality of $l$.

For the other direction, assume that each pair $g,h\in G$ has a maximal common divisor. Let $H:=\{h_1,\ldots,h_n\}\subseteq G$ be a finite non-empty set, and let $d$ be an upper bound of $H$. We now have
\[
d=g_1h_1=\cdots=g_nh_n
\]
for some $g_1,\ldots,g_n\in G$. Let $g'\in G$ be a maximal common divisor of $g_1,\ldots,g_n$. Hence there exists $g_i'\in G$ such that $g_i=g_i'g'$ for all $i\in\{1,\ldots,n\}$. Also, $d=d'g'$ for some $d'\in G$. It is now easy to see that the maximal common divisor of $g'_1,\ldots,g'_n$ is $1$, and that $d'$ is a minimal upper bound of $H$.
\end{proof}

\begin{notation}\label{hattu}
Let $S\subseteq \C$ be a finite subset. We denote the set of minimal upper bounds of non-empty subsets of $S$ by
\[
\hat{S} := \bigcup_{\emptyset\neq S'\subseteq S} \mub_\C(S').
\]
We notice that if $\C$ is weakly bounded from above, then $\hat{S}$ is finite.
\end{notation}

Using the terminology from~\cite{Corbet}, a set $S\subseteq \C$ is \emph{a framing set of $M$} if every $c\in{\uparrow}\supp(M)$ has an element $s\in S\cap{\downarrow} c$, called a \emph{frame of $c$}, such that $M(s\leq c')$ is an isomorphism for all $s\leq c'\leq c$.

\begin{lemma}\label{tahti}
If an $R\C$-module $M$ has a framing set $S$, then $M$ is $S$-determined.

Conversely, if $\C$ is weakly bounded from above and mub-complete, and $M$ is $S$-determined for some finite set $S\subseteq \C$, then $\hat{S}$ is a finite framing set of $M$. In particular, if $c\in \C$, then every $s\in\mub(S \cap{\downarrow} c)\subseteq \hat S$ is a frame of $c$ such that $S\cap {\downarrow} c =  S \cap {\downarrow} s$.
\end{lemma} 

\begin{proof}
Assume first that $S$ is a framing set for $M$. If $c\in \supp(M)$, then there exists a frame $s\in S$ of $c$, and therefore $c\in {\uparrow} S$. Thus $\supp(M)\subseteq {\uparrow} S$. Let $c\leq d$ in $\C$ such that $S\cap{\downarrow} c=S\cap {\downarrow} d$. If $d\notin {\uparrow} \supp(M)$, we have $M(c)=M(d)=0$, and we are done. Otherwise, there exists a frame $s\in S$ of $d$. Since $S\cap{\downarrow} c=S\cap {\downarrow} d$, we see that $s\leq c\leq d$. Therefore $M(c\leq d)$ is an isomorphism.

Assume next that $\C$ is weakly bounded from above and mub-complete, and $M$ is $S$-determined for some finite set $S$. Since $\C$ is weakly bounded from above, $\hat{S}$ is finite. Let $c\in{\uparrow} \supp(M)$. Now there exists an element $b\leq c$ such that $M(b)\neq 0$. Since $M$ is $S$-determined, we see that $S\cap {\downarrow} c \supseteq S\cap{\downarrow} b\neq \emptyset$. Thus $c$ is an upper bound of the non-empty set $S\cap{\downarrow} c$, so by mub-completeness there exists a minimal upper bound $s\in \mub(S\cap {\downarrow} c)\subseteq \hat{S}$ such that $s\leq c$. It follows that $S\cap {\downarrow} c \subseteq S \cap {\downarrow} s$. Obviously $s\leq c$ implies $S\cap{\downarrow} s \subseteq S\cap{\downarrow} c$. Hence $S\cap{\downarrow} s = S\cap{\downarrow} c$. If $s\leq c'\leq c$, then trivially ${\downarrow} S\cap s={\downarrow} S\cap c'$, so $M(s\leq c')$ is an isomorphism. 
\end{proof}


\subsection{Finitely presented \texorpdfstring{$R\C$}{RC}-modules in mub-complete posets}

In the next proposition we find out how the minimal upper bounds connect to births and deaths relative to $S$. This allows us to prove our main result, Theorem~\ref{fin_pres}.

\begin{proposition}\label{syntyma}
Let $\C$ be weakly bounded from above and mub-complete. Let $M$ be an $R\C$-module that is $S$-determined for some finite $S\subseteq \C$. If $B_{\hat{S}}(M)\subseteq S$, then $D_{\hat{S}}(M)\subseteq \hat{S}$.
\end{proposition}

\begin{proof}
Suppose that $B_{\hat{S}}(M)\subseteq S$. This implies that $M$ is $\hat{S}$-generated by Proposition~\ref{FSP_apu}, so $M$ is $B_{\hat{S}}(M)$-generated by Proposition~\ref{syntyma_minimi}. Let $c\in D_{\hat{S}}(M)$, so that $\lambda_{M,c}$ is not a monomorphism.  This means that there exist $d_1,\ldots,d_n\in \hat{S}$ such that $d_i < c$ for all $i=\{1,\ldots,n\}$, and non-zero elements $x_1\in M(d_1),\ldots,x_n\in M(d_n)$ for which $\sum_{i=1}^n M(d_i\leq c)(x_i)=0$. If there exists $c'\in \hat{S}$ such that $d_i\leq c' <c$ for all $i\in \{1,\ldots,n\}$, we may assume that $\sum_{i=1}^n M(d_i\leq c')(x_i)\neq 0$, since the homomorphism
\[
\bigoplus_{i=1}^n M(d_i)\rightarrow \colim_{d< c, \ d\in \hat{S}} M(d)
\]
factors through $M(c')$. Because $M$ is $B_{\hat{S}}(M)$-generated, we may also assume that $d_i\in B_{\hat{S}}(M)\subseteq S$ for all $i\in \{1,\ldots,n\}$.

On the other hand, by Lemma~\ref{tahti}, we have $S\cap {\downarrow} s = S \cap {\downarrow} c$ for some frame $s\in \hat{S}$ of $c$. This implies that $d_i\leq s$ for all $i\in \{1,\ldots,n\}$. If $s<c$, we get a contradiction $\sum_{i=1}^n M(d_i\leq s)(x_i)=0$. Therefore $c=s\in \hat{S}$.
\end{proof}

\begin{corollary}\label{tuplahattu}
Let $\C$ be weakly bounded from above and mub-complete. Let $M$ be an $S$-determined $R\C$-module, where $S \subseteq \C$ is a finite subset.
Then $\hat{\hat S}$ is an FSP of $M$.
\end{corollary}

\begin{proof}
Obviously $M$ is $\hat{S}$-determined, since $S\subseteq \hat{S}$. By Proposition~\ref{syntyma}, it is enough to show that $B_{\hat{\hat S}}(M)\subseteq \hat{S}$, because then $D_{\hat{\hat S}}(M)\subseteq \hat{\hat S}$. The rest now follows from Proposition~\ref{FSP_apu}.
Since $\hat{S}\subseteq \hat{\hat S}$, by Remark~\ref{syntyma_inkluusio} we have $B_{\hat{\hat S}}(M)\subseteq B_{\hat S}(M)$. Let $c\in \C$. If $c\notin \hat{S}$, then $c$ has a frame $s\in \hat{S}$ by Lemma~\ref{tahti}. This means that $c\notin B_{\hat{S}}(M)$, and thus $B_{\hat S}(M)\subseteq \hat{S}$.
\end{proof}

We sum up Proposition~\ref{FSP_ehdot} and Corollary~\ref{tuplahattu} in the next corollary.

\begin{corollary}\label{FSP_luonnehdinta}
Let $\C$ be weakly bounded from above and mub-complete. An $R\C$-module $M$ has an FSP if and only if $M$ is $S$-determined for some finite $S\subseteq \C$.
\end{corollary}

Finally, we get our new characterization of finitely presented modules.

\begin{theorem}\label{fin_pres}
Let $M$ be an $R\C$-module. If $M$ is finitely presented, then
\begin{itemize}
\item[1)] $M(c)$ is finitely presented for all $c\in \C$
\item[2)] $M$ is $S$-determined for some finite $S\subseteq \C$.
\end{itemize}
Furthermore, if $\C$ is weakly bounded from above and mub-complete, and $M$ satisfies conditions 1) and 2), then $M$ is finitely presented.
\end{theorem}

\begin{proof}
Using Proposition~\ref{roiter_theorem}, the first part immediately follows from Proposition~\ref{FSP_ehdot} and the second part from Corollary~\ref{tuplahattu}.
\end{proof}

\begin{remark}\label{scola_corbet}
Let $G$ be a monoid. 
Theorem~\ref{fin_pres} and Lemma~\ref{tahti} show us that the $RG$-modules of finitely presented type of Corbet and Kerber~(\cite{Corbet}*{p.~19, Def.~15}) are the same thing as finitely presented $RG$-modules. 

Furthermore, let $A$ be a free $G$-act that is mub-complete and weakly bounded from above as a poset. Starting from the isomorphism $RA\text{-}\textbf{Mod} \cong A\text{-gr } R[G]\text{-}\textbf{Mod}$ of Corollary \ref{lineaari ekvivalenssi}, and using the fact that being finitely presented is a categorical property, we get an isomorphism between finitely presented $RA$-modules and finitely presented $A$-graded $R[G]$-modules. Taking $A=G$ now gives the commutative case of \cite{Corbet}*{p.~25, Thm.~21}.  
\end{remark}

\begin{example}\label{mub_esimerkki}
Let $\C=\{a,b\}\cup \mathbb Z$, where $a< n$ and $b<n$ for all $n\in \mathbb Z$, and $\mathbb Z$ has the usual ordering. Then the $R\C$-module $M:=R[\Mor_\C(a,-)]\oplus R[\Mor_\C(b,-)]$ is obviously finitely presented, but $M$ does not have a finite framing set even though it is $\{a,b\}$-determined. Caution is required here: If we define an $R\C$-module $N$ by
\[
N(a)=N(b)=R \quad \text{and} \quad N(n)=R^3,
\]
for all $n\in \mathbb Z$, then $N$ satisfies the conditions $1)$ and $2)$ in Theorem~\ref{fin_pres} but is not finitely presented. This follows from the fact that $\C$ is not mub-complete.
\end{example}



\end{document}